\documentclass [winedt,yap]{iitparc}
\usepackage{cite}
\usepackage{amsmath,amssymb,amsfonts,bm}
\usepackage{eufrak}
\usepackage{lscape}
\usepackage{floatfig,wrapfig,epsfig}
\usepackage{subfigure}
\usepackage{color}
\usepackage{psboxit}
\usepackage{rotating}
\usepackage{curves}
\usepackage{ulem}
\usepackage[mathscr]{eucal}
\usepackage{array}            
\RequirePackage{psfrag}
\RequirePackage{graphicx}

\begin{document}\normalem
\initfloatingfigs
\frontmatter          

\IssuePrice{25.00}%
\TransYearOfIssue{2017}%
\TransCopyrightYear{2017}%
\OrigYearOfIssue{2017}%
\OrigCopyrightYear{2017}%

\TransVolumeNo{78}%
\TransIssueNo{1}%
\OrigIssueNo{1}%
\OrigPages{106--120}


\mainmatter

\setcounter{page}{92}
\CRubrika{INTELLECTUAL CONTROL SYSTEMS}
\Rubrika{INTELLECTUAL CONTROL SYSTEMS}

%
%
%
\def\rank{\operatorname{rank}}                   
\def\diag{\operatorname{diag}}                   
\def\bp{\bm\pi}
%
%

\Rubrika{Intellectual control systems}

\title{Models of latent consensus%
\thanks{The work of the second author was supported by the Russian Science Foundation (project
16-11-00063 granted to IRE RAS).}}

\author{R. P. Agaev and        
        P. Yu. Chebotarev}     
\authorrunning{Agaev, Chebotarev}
\OrigCopyrightedAuthors{R.P. Agaev and P.Yu. Chebotarev}
\institute{Trapeznikov Institute of Control Sciences, Russian Academy of Sciences, Moscow, Russia}

\received{Received January 16, 2016}

\maketitle

\begin{abstract}
The paper studies the problem of achieving consensus in multi-agent systems in the case where the dependency digraph~$\Gamma$ has no spanning in-tree. We consider the regularization protocol that amounts to the addition of a dummy agent (\emph{hub}) uniformly connected to the agents. The presence of such a hub guarantees the achievement of an asymptotic consensus. For the ``evaporation'' of the dummy agent, the strength of its influences on the other agents vanishes, which leads to the concept of \emph{latent consensus}. We obtain a closed-form expression for the consensus when the connections of the hub are symmetric; in this case, the impact of the hub upon the consensus remains fixed. On the other hand, if the hub is essentially influenced by the agents, whereas its influence on them tends to zero, then the consensus is expressed by the scalar product of the vector of column means of the Laplacian eigenprojection of $\Gamma $ and the initial state vector of the system. Another protocol, which assumes the presence of vanishingly weak uniform background links between the agents, leads to the same latent consensus.

\end{abstract}

\setcounter{footnote}{1}

\section{Introduction}

In networked multi-agent systems (MAS) consensus is asymptotically reached for any initial state if and only if $0$ is a simple eigenvalue of the Laplacian matrix of the weighted dependency digraph. An equivalent condition is the presence of a spanning in-tree in this digraph~\cite{Olfati-SaberMurray04, MesbahiEgerstedt10book}. Suppose that these conditions are violated. Is it still possible to associate with the system a certain ``most natural'' consensus achievable under a minimal or even infinitesimal transformation of the system?

In this paper, we propose two approaches to the determination of such a consensus, which can be referred to as a \emph{latent consensus}. At the first approach, the system is supplemented by a dummy agent called a \emph{hub} with a vanishingly weak uniform influence on the agents, which guarantees the achievement of an asymptotic consensus. At the second approach, the system is supplemented by uniform and vanishingly weak ``background'' links between the existing agents. It turns out that both approaches lead to the same latent consensus. If the additional links are not vanishingly weak, then in both cases, consensus is the average of the initial states (opinions) with weights equal to the column means of the parametric matrix of spanning in-forests of the dependency digraph. When the strength of the additional links tends to zero, the weights of the agents' states converge to the column means of the matrix of \emph{maximum} in-forests. The closed-form expressions for the latent consensus in both continuous and discrete networked multi-agent systems involve the eigenprojection of the Laplacian matrix of the dependency digraph.

A basically different approach to solving the same problem is achieving consensus not due to a small change in the structure of relations, but due to a minimal appropriate correction of the initial state. This approach is realized by the \emph{method of orthogonal projection\/} on the consensus subspace~\cite{AgaChe11ARC1,AgaChe15ARC}. Its results are generally different from those of the methods of \emph{eigenprojection\/} presented in this paper.

The concept of latent consensus requires a discussion. The consensus problem admits two interpretations.
In the first one, we synthesize a coordination protocol for a technical multi-agent system, e.g., averaging the speeds of drones flying in a flock. In the second interpretation, we model smooth or iterative convergence of opinions of the respondents (or experts) due to their influences on each other. In the first interpretation, there are no unreported connections between agents: just those included in the protocol matter. In contrast, the second interpretation allows for the presence of background, latent ``public opinion'' affecting agents along with their interactions specified explicitly. This ``public opinion'' can be simulated in two ways: (1) by  rather weak uniform pairwise connections between agents and (2) by introducing a special ``invisible'' agent uniformly affecting all real agents and experiencing, in turn, their influence. Consensus achievable in these models in the case where the explicit links do not guarantee its existence can be called a \emph{latent consensuses\/}, which justifies the use of this term.

Returning to the first interpretation, i.e., to a technical system, the introduction of a dummy hub can be seen as the assumption of the existence of a ``center'' with weak uniform influences on the agents. An alternative, i.e., adding weak background pairwise interactions, can hardy be recommended as a tool to ensure consensus, as it would dramatically increase the computational complexity for systems with large and relatively sparse matrices of influences. On the other hand, in virtue of the results of this paper, the latter approach does not give anything new, since it leads to the same consensus as the addition of a dummy hub.

Stabilization by the addition of weak interactions can be treated as an embodiment of the idea of low gain feedback~\cite{Lin99}. This approach has been used as applied to consensus problems, e.g., for system with general linear dynamics~\cite{SeoShimBack09,SuChenLamLin13} and with delay~\cite{Olfati-SaberMurray04}. As distinct from these works, which studied the means to achieve consensus, in the present case of the simplest dynamics, asymptotic consensus is guaranteed by the presence of a spanning in-tree in the dependency digraph, while the object of interest is the functional relationship of the consensus with the invariants of the digraph. The latter is the subject of this paper.

The paper is organized as follows. 
In sections~\ref{s_model} and~\ref{hub_bezv}, we obtain expressions for the consensus in systems with a dummy agent (``hub'') that weakly influences the existing agents and, in turn, is experiencing their impact. In Section~\ref{s_backgr}, we study protocols with weak ``background'' links between agents. Section~\ref{s_ortho} contains an example used to compare the methods of regularization based on the addition of vanishingly weak links with the method of orthogonal projection. Finally, in Section~\ref{s_discr} we obtain closed-form expressions for the consensus in weak links regularization protocols applied to DeGroot's iterative pooling process.

\section{Notation and auxiliary results}
\label{s_notat090213}

Consider the continuous protocol of consensus seeking in a multi-agent system \cite{Olfati-SaberMurray04,RenCao11distributed}
\begin{gather}\label{e_model0715}
\dot{x}_i(t)=-\sum_{j=1}^na_{ij}\left({x_i(t)-x_j(t)}\right),\quad
i=1,\ldots,n ,
\end{gather}
where $x_i(t)$ is the characteristic to be adjusted of agent~$i.$

With the \emph{dependency matrix\/} $A=(a_{ij})_{n\times n}$ of protocol \eqref{e_model0715} we associate the \emph{weighted dependency digraph} $\Gamma $ with vertex set $\{1,\ldots,n\}$ that for every element $a_{ij}>0$ of $A$ has an arc $(i,j)$ with weight $w_{ij}=a_{ij}$. Thus, arcs in $\Gamma $ are drawn in the direction of dependency: from each dependent vertex / agent $i$ to every vertex / agent $j$ affecting~$i$; the weight $w_{ij}$ of arc $(i,j)$ reflects the strength of this influence.
The \emph{Laplacian matrix\/} of the weighted digraph $\Gamma $ is defined \cite{CheAga02ap} by 
\begin{gather}\label{e_modeL}
L=\diag(A\bm1)-A,
\end{gather}
where $\bm1=(\underbrace{1,\ldots, 1}_n)^{\mathrm T},$ $\diag(\boldsymbol{v})$ being the operator transforming a vector $\boldsymbol{v}$ to the diagonal matrix with $\boldsymbol{v}$ on the diagonal. The class of Laplacian matrices of weighted digraphs coincides with the class of square matrices with nonpositive off-diagonal entries and zero row sums. In what follows, the term \emph{Laplacian matrix\/} will refer to an element of this class.

It is easy to verify that matrix \eqref{e_modeL} appears in the matrix representation of the model~\eqref{e_model0715}:
\begin{gather}\label{e_040715q}
\dot {\boldsymbol{x}}(t)=-L\,\boldsymbol{x}(t),
\end{gather}
 where $\boldsymbol{x}(t)=(x_1(t),\ldots, x_n(t))^{\mathrm T}.$

Now recall the necessary notation of matrix theory. Suppose that ${A\in\mathbb{C}^{n\times n}}$ is any square matrix,
${\mathcal{R}}({A})$ and ${\mathcal{N}}({A})$ being the range and the null space of~$A,$ respectively. Let ${\nu={\mathrm{ind}\,} A}$ be the {\it index\/} of $A,$ i.e., the smallest ${k\in\{0,1,\ldots\}}$ such that ${\rank A^{k+1}=\rank A^k}$ (by definition, we assume ${A^0\equiv I}$, where $I$ is the identity matrix of order~$n$).

The {\it eigenprojection\/} of matrix $A,$ {\it corresponding to the eigenvalue}~$0$ (or simply the {\it eigenprojection of\/}~$A$) is a projection (i.e., an idempotent matrix) $A^\vdash$ such that
${\mathcal{R}}(A^\vdash)={\mathcal{N}}(A^\nu)$ and
${\mathcal{N}}(A^\vdash)={\mathcal{R}}(A^\nu).$ Such a matrix $A^\vdash$ is always unique.

\begin{proposition}\label{p_limInf}
If $\boldsymbol{x}(t)$ is a solution to the system of equations \eqref{e_040715q}$,$ then
$$
\lim_{t\to\infty}\boldsymbol{x}(t)=L^\vdash  \boldsymbol{x}(0).
$$
\end{proposition}

Proposition~\ref{p_limInf} follows from the forest consensus theorem (Theorem~1 in~\cite{CheAga14IEEETAC}).

\begin{lemma}[{\hspace*{-1mm}\cite{CheAga02ap}}] 
\label{l_Mforest} The eigenprojection $L^\vdash$ coincides with the normalized matrix ${\bar{J}=(\bar{J}_{ks})}$ of
\emph{maximum in-forests\/} of the dependency digraph $\Gamma $\/$:$
$$
\bar{J}_{ks}=\frac{f_{ks}}f,\quad k,s=1,\ldots,n ,
$$
where $f$ is the total weight\/\footnote{The weight of a digraph (e.g., of an in-forest) is the product of the weights of all its arcs.} of all maximum in-forests of digraph $\Gamma ,\,$ $f_{ks}$ being the total weight of those of them that have $k$ belonging to a tree with root $($sink$!)$~$s.$                                            
\end{lemma}

In the subsequent proofs, we make use of the following lemma.
\begin{lemma}[\cite{AgaChe00'}]\label{02072015th1}
$L^\vdash=\lim_{\tau\to\infty}(I+\tau L)^{-1}.$
\end{lemma}

Now we introduce some additional notation\footnote{Square brackets will denote row vectors and row representations of block matrices whose entries (blocks) are separated with spaces rather than commas.}:
$\bm0=(\underbrace{0,\ldots, 0}_n)^{\mathrm T},$
$E=\frac1n\bm1\bm1^{\mathrm T},$ $\bm1'=[\bm1^{\mathrm T}\;\,1]^{\mathrm T},$ $I'$ and $E'$ are two matrices of order $n+1$:
the identity matrix and matrix with entries $\tfrac1{n+1},$ respectively;
\[L_0=\!\begin{pmatrix}
      L &\bm0\\[.3em]
\bm0^{\mathrm T} & 0    \end{pmatrix},\quad
I_0=\!\begin{pmatrix}
      I &\bm0\\[.3em]
\bm0^{\mathrm T} & 0     \end{pmatrix}.
\]

The proofs of the main results of this paper rely on the following lemma.
\begin{lemma}\label{31102015pr1}
If for matrices $A$ and $C,$ ${AC=A}$ and ${C^2=C}$ hold, ${\delta>0},$ and the real parts of the eigenvalues of ${\delta I+A}$ are positive\/$,$ then
$$
\lim_{t\to\infty}e^{-\left(A+\delta C\right)t} = (I-C)\big(I+\delta^{-1}  A\big)^{-1}. 
$$
\end{lemma}

The proof of Lemma~\ref{31102015pr1} and the other statements are given in the Appendix.

We will also use the following result. 
\begin{lemma}\label{l_sopL}
If for Laplacian matrices $A$ and $C,$ ${AC=A}$ and ${C^2=C}$ hold, then 
${(A+\delta C)^\vdash =}\linebreak{=(I-C)\big(I+\delta^{-1} A\big)^{-1}\;}$ for any~${\delta>0}$.
\end{lemma}

The proof of this lemma is based on Lemma~\ref{31102015pr1} and Lemma~\ref{l_sopLe}.
\begin{lemma}\label{l_sopLe}
For any Laplacian matrix $L,$ it holds that ${L^\vdash =\lim_{t\to\infty}e^{-L  t}}$.
\end{lemma}

\section{Consensus protocols with a symmetric hub}
\label{s_model}

Consider the following protocols of consensus seeking:
\begin{gather}\label{e_modeLK040715c} \dot{
\boldsymbol{y}}(t)=
-(L_0+H_{\delta,\boldsymbol{v}})\,\boldsymbol{y}(t),
\end{gather} 
where $\boldsymbol{y}(t)=(y_1(t),\ldots, y_{n+1}(t))^{\mathrm T}$, $H_{\delta  ,\boldsymbol{v}}=\delta   H_I+H_{\boldsymbol{v}}$,
$\delta  >0$,
\begin{gather}\label{e_W1} H_I =
\begin{pmatrix}
I       &-\bm1\\[.3em]
\bm0^{\mathrm T} &    0   \end{pmatrix}\quad \mbox{and}\quad
H_{\boldsymbol{v}} =
\begin{pmatrix}
0_{n\times n} & \bm0\\[.3em]
   -\boldsymbol{v}^{\mathrm T}     &   s \end{pmatrix}
\end{gather}
are Laplacian matrices of order $n+1,\,$ $\boldsymbol{v}=(v_1,\ldots, v_n)^{\mathrm T},\,$ and $\,s=\sum_{i=1}^n v_i.$

Every protocol \eqref{e_modeLK040715c} is obtained from \eqref{e_040715q} by adding an ${(n+1)}$st agent that influences the $n$ agents of the system with intensity~$\delta  $ and is subjected to their influences whose intensities are given by the vector~$\boldsymbol{v}.$ This supplementary agent (and the corresponding vertex of the dependency digraph) will be called a \emph{hub}\/\footnote{In \cite{CheAga02ap}, before Definition~3, we gave references to a number of publications that use this concept either explicitly or implicitly.}.
The motivation for the introduction of a hub is to assure the presence of a spanning in-tree in the dependency digraph (which guarantees asymptotic consensus) due to hub's uniform weak influences on the agents. 

The eigenprojection $(L_0+H_{\delta  ,\boldsymbol{v}})^\vdash$ has a simple representation which provides an expression for the consensus by means of Proposition~\ref{p_limInf}.

\begin{lemma}\label{02072015b}
$\,(L_0+H_{\delta  ,\boldsymbol{v}})^\vdash
=\bm1'\tfrac1{s+\delta  } \Big[\boldsymbol{v}^{\mathrm T}\big(I+\tfrac1\delta   L\big)^{-1}\;\;\, \delta  \Big].$
\end{lemma}

\begin{theorem}\label{t_UniHub}
In the system \eqref{e_modeLK040715c} with a hub that uniformly influences the agents$,$ the consensus is $\:\lim_{t\to\infty}\boldsymbol{y}(t) 
=\bm1'\tfrac1{s+\delta  } \Big[\boldsymbol{v}^{\mathrm T}\big(I+\tfrac1\delta   L\big)^{-1}\;\; \delta \Big]\boldsymbol{y}(0).$
\end{theorem}

Now consider a hub whose links are symmetric: ${v_i=\delta  }$,
${i=1,\ldots, n}$. Then $\boldsymbol{v}=\delta  \bm1$ and $s=\delta n.$ Using Lemmas~\ref{02072015b} and \ref{02072015th1}, we obtain:
{\arraycolsep=.4mm\begin{eqnarray}\label{e_symmHub} (L_0+H_{\delta
,\,\delta \bm1})^\vdash
&=&\bm1'\frac1{n+1} \left[\bm1^{\mathrm T}\left(I+\dfrac1\delta   L\right)^{-1}\;\; 1\right],\\[.3em]
\lim_{\delta  \to+0}(L_0+H_{\delta  ,\,\delta  \bm1})^\vdash
&=&\bm1'\frac1{n+1} \Big[\bm1^{\mathrm T}\!\bar{J}\;\;\, 1\Big].
\label{e_symmHublim} \end{eqnarray} }%

Eqs.\ \eqref{e_symmHub} and \eqref{e_symmHublim} and Proposition~\ref{p_limInf} imply Proposition~\ref{p_limInf+}.

\begin{proposition}\label{p_limInf+}
If $\boldsymbol{y}(t)$ satisfies \eqref{e_modeLK040715c} and $\boldsymbol{v}=\delta  \bm1,$ then
{\arraycolsep=.4mm\begin{eqnarray*}
\lim_{t\to\infty}\boldsymbol{y}(t)
&=& \bm1'\dfrac1{n+1}\! \left(\sum_{j=1}^ns_j^\delta  \, y_j(0)+y_{n+1}(0)\right),\\
\lim_{\delta  \to+0}\lim_{t\to\infty}\boldsymbol{y}(t)
&=&\bm1'\dfrac1{n+1}\! \left(\sum_{j=1}^n\bar{J}_{ \cdot
j}y_j(0)+y_{n+1}(0)\right),
\end{eqnarray*}}
where $s_j^\delta  \,$ and $\,\bar{J}_{\cdot j}=\sum_{i=1}^n\bar{J}_{ij}$ are the $j$th column sums of matrices
$\,(I+\tfrac1\delta   L\big)^{-1}$ and $\,L^\vdash=\bar{J}=(\bar{J}_{ij})_{n\times n},\,$ respectively.
\end{proposition}

By Proposition~\ref{p_limInf+}, the weight of the initial state $y_{n+1}(0)$ of a hub with symmetric links is $(n+1)^{-1},$ irrespective of~$\delta  .$

\begin{corollary}\label{co_limlim}
$1.$ If the system \eqref{e_040715q} has an asymptotic consensus expressed by the scalar product of the vectors
$\boldsymbol{x}(0)$ and $(v_1,\ldots, v_n)^{\mathrm T},$ then the asymptotic consensus in the system \eqref{e_modeLK040715c} with
${\boldsymbol{v}=\delta  \bm1}$ and~${\delta  \to+0}$ is the scalar product of the vectors $\boldsymbol{y}(0)$ and
${(n+1)^{-1}(nv_1,\ldots, nv_n,1)^{\mathrm T}}$.

$2.$ If the initial state of the hub satisfies the equality $\,{y_{n+1}(0)=\frac1n\big(\sum_{i,j=1}^n\bar{J}_{ij}y_j(0)\big)},$ then this equality also expresses the asymptotic consensus in the protocol~\eqref{e_modeLK040715c} with ${\boldsymbol{v}=\delta\bm1}$ and ${\delta  \to+0}.$

$3.$ The coincidence of the asymptotic consensuses of systems \eqref{e_040715q} and \eqref{e_modeLK040715c} in item\;$1$ takes place if and only if $y_{n+1}(0)$ satisfies the equality in item\;$2.$
\end{corollary}

The proof of Corollary~\ref{co_limlim} is obvious.

According to Proposition~\ref{p_limInf} and Lemma~\ref{l_Mforest}, the consensus in protocols~\eqref{e_modeLK040715c}
is the weighted average of the initial states of the agents (including the hub), where the weight of the state of agent $j$
is the specific weight of spanning in-trees with sink $j$ in the dependency digraph whose Laplacian matrix is
${L_0+H_{\delta  ,\boldsymbol{v}}.}$ By Theorem~\ref{t_UniHub} and the matrix forest theorem \cite{CheSha97'} 
this weight can also be found by averaging (with weights specified by~$\boldsymbol{v}$) the specific $(1/\delta )$-weights\footnote{By
``$b$-weight'' of a weighted digraph we mean the product of its arc weights, each multiplied by~$b.$} of the forests, converging to $j$ from all the vertices of~$\Gamma .$

Let the links of the hub be symmetric, i.e., ${\boldsymbol{v}=\delta \bm1}$.
According to Proposition~\ref{p_limInf+}, for any $\Gamma $ and ${\delta  \to\infty},$ the consensus equals ${\frac1{n+1}\sum_{j=1}^{n+1} y_j(0)},$ the arithmetic mean of the initial states of all agents, including the hub.
As ${\delta  \to+0},$ to determine the weight of the initial state of agent ${j\le n},$ the average specific weight
$\,\bar{J}_{\cdot j}/n$ of the \emph{maximum\/} forests converging to vertex $j$ in $\Gamma $ must be taken with the coefficient $\tfrac n{n+1},$ and the initial state of the hub with weight $(n+1)^{-1}.$ Such a consensus cannot be termed ``latent'',
as it essentially depends on the initial state of the hub added to the system, which is an \emph{external\/} parameter.

\section{A hub subordinate to the agents}
\label{hub_bezv}

According to Proposition~\ref{p_limInf+}, a symmetric hub in protocols \eqref{e_modeLK040715c} outwardly affects the consensus with an average strength: the weight of its initial state is $(n+1)^{-1}.$ The constancy of this weight, even when the strength of hub's links goes to zero, may seem strange. The reason is as follows. If a symmetric hub very weakly depends on the agents, then the influence of its initial state is though weak, but rather stable. So the resulting coefficient $(n+1)^{-1}$ of its influence can be treated as a manifestation of the principle of ``little strokes fell great oaks''. 

Therefore, if our aim is to ensure consensus \emph{without arbitrary external influences\/}, then the links of the hub should be 
set in a different way. Consider the protocols
\eqref{e_modeLK040715c}, where
\begin{gather}\label{e_trueHub} \delta
/s\to+0\;\mbox{ as }\;\delta  \to+0.
\end{gather}

Thereby, we consider the process, in which the influence $\delta $ of the hub on the agents becomes vanishingly small not only in itself, but also compared to \emph{their\/} influence $s$ on the hub. This ``makes the hub forget'' its initial state, so that its impact on the agents becomes a translation of their own average influence (a kind of ``public opinion'').

A counterpart of Eq.~\eqref{e_symmHublim} for this case is the following statement.

\begin{lemma}\label{l_SopNs}
Under \eqref{e_trueHub} and the existence of ${\tilde{\boldsymbol{v}}=\lim_{\delta \to+0}\tfrac1s\boldsymbol{v}},$ it holds that
$
{\lim_{\delta  \to+0}(L_0+H_{\delta  ,\,\boldsymbol{v}})^\vdash}\linebreak{
=\bm1'\big[\tilde{\boldsymbol{v}}^{\mathrm T}\!\bar{J}\;\;\, 0\big]}
$.
\end{lemma}

Along with Lemma~\ref{l_SopNs}, Proposition~\ref{p_limInf} implies Theorem~\ref{t_depsop}. 
\begin{theorem}\label{t_depsop}
Under the conditions of Lemma~$\ref{l_SopNs},$ for protocols~\eqref{e_modeLK040715c} we have\/$\/:$

$1.$ $\lim\limits_{\delta
\to+0}\lim\limits_{t\to\infty}\boldsymbol{y}(t) =
\bm1'\big[\tilde{\boldsymbol{v}}^{\mathrm T}\!\bar{J}\;\; 0\big]
\boldsymbol{y}(0)$$;$

$2.$ If ${\boldsymbol{v}=\tfrac s n\bm1},$ then
${\lim\limits_{\delta\to+0}\lim\limits_{t\to\infty}\boldsymbol{y}(t) 
=\bm1'\big[\tfrac1n\bm1^{\mathrm T}\!\bar{J}\;\:0\big]\boldsymbol{y}(0) 
=\bm1'\,\tfrac1n\sum_{j=1}^n\bar{J}_{ \cdot j}y_j(0)},$ 
where ${\bar{J}_{\cdot j}=\sum_{i=1}^n\bar{J}_{ij}.}$
\end{theorem}

Consensus in item\:2 of Theorem~\ref{t_depsop} can legitimately be termed \emph{latent}\/: it takes place for a wide class of
protocols \eqref{e_modeLK040715c} with a hub subordinate (in accordance with \eqref{e_trueHub}) to the agents and uniformly, but infinitesimally affecting them. This consensus equals the weighted average of the initial states of $n$ agents with weights equal to the  column means of $\bar{J}=L^\vdash.$
Herewith, due to Lemma~\ref{l_Mforest} the weight of the initial state of agent $j$ is the average (over $i$) specific
weight of the maximum in-forests in $\Gamma $, where $j$ is a sink accessible from~$i.$ It is easily seen that this weight equals the probability to access sink $j$ from a randomly chosen vertex in the model of motion along maximum in-forests chosen with probabilities
proportional to their weights.

\section{Protocols with weak background links between the agents}
\label{s_backgr}

Now we turn to a different method of regularization. Consider the protocols of consensus seeking in which the initial dependency digraph $\Gamma $ is supplemented by the complete digraph of ``background links'' with uniform arc weights:
\begin{gather}\label{e_modeLK011115c} 
\dot{ \boldsymbol{x}}(t)
=-(L+\delta K)\,\boldsymbol{x}(t), 
\end{gather}
where $\delta>0$ and $K=I-E$ ($K$ is the Laplacian matrix of the complete digraph with arc weights~$1/n$). 

First, in the same way as before, we study a more general type of regularization.
Let $\boldsymbol{v}$ be a distribution vector on the set of agents, i.e., ${v_i\ge 0}$ (${i=1,\ldots, n}$) and
${\sum_{i=1}^nv_i=1}$. Set 
\begin{gather}\label{e_modeL-PR}
\dot{ \boldsymbol{x}}(t)=-(L+\delta   D)\,\boldsymbol{x}(t),
\end{gather} 
where $\,\delta>0,$ $\,D=I-V,$ and $\,V=\bm1\boldsymbol{v}^{\mathrm T}.$ $D$ is Laplacian; $\delta
v_i$ is the influence of agent $i$ on any other.

\begin{lemma}\label{l_PRsop}
In the above notation\/$,$

$1.$ $ (L+\delta   D)^\vdash = \bm1 \boldsymbol{v}^{\mathrm T}
\big(I+\tfrac1\delta   L\big)^{-1}$$;$

$2.$ $ \lim_{\delta  \to0}(L+\delta   D)^\vdash =\bm1
\boldsymbol{v}^{\mathrm T}\!\bar{J}.$
\end{lemma}

Proposition~\ref{p_limInf} and Lemma~\ref{l_PRsop} imply Theorem~\ref{t_limInf_0111a}.

\begin{theorem}\label{t_limInf_0111a}
If $\boldsymbol{x}(t)$ is a solution to the system of equations \eqref{e_modeL-PR}$,$ then 

$1.$ $\lim_{t\to\infty}\boldsymbol{x}(t) \!=\! \bm1
\boldsymbol{v}^{\mathrm T} \big(I+\tfrac1\delta
L\big)^{-1}\boldsymbol{x}(0)$$;$

$2.$ $\lim_{\delta  \to+0}\lim_{t\to\infty}\boldsymbol{x}(t) =\bm1
\boldsymbol{v}^{\mathrm T}\!\bar{J} \boldsymbol{x}(0).$
\end{theorem}

Observe that with ${\boldsymbol{v}=\tfrac1n\bm1}$ we have ${\bm1\boldsymbol{v}^{\mathrm T}=E}$ and ${D=K}$, and so Corollary~\ref{c_limInf} holds true.
\begin{corollary}\label{c_limInf}
If $\boldsymbol{x}(t)$ is the solution to the system of equations
\eqref{e_modeLK011115c}$,$ then

$1.$ $\lim_{t\to\infty}\boldsymbol{x}(t)
\!=\!\bm1\,\tfrac1n \sum_{j=1}^ns_j^\delta    x_j(0),$ where
$s_j^\delta  $ are the column sums of $(I+\tfrac1\delta L\big)^{-1}$$;$

$2.$ $\lim_{\delta  \to+0}\lim_{t\to\infty}\boldsymbol{x}(t)
=\bm1\,\tfrac1n \sum_{j=1}^n\bar{J}_{ \cdot j}x_j(0),$ where
$\bar{J}_{\cdot j}=\sum_{i=1}^n\bar{J}_{ij},\; j=1,\ldots, n.$
\end{corollary}

Thus, in the protocols with additional uniform background links~\eqref{e_modeLK011115c}, consensus is determined by averaging the column entries of $(I+\tfrac1\delta L\big)^{-1},$ which is the parametric matrix of in-forests \cite{CheSha97'} of digraph~$\Gamma .$ Furthermore, if the arc weights of the added complete digraph tend to zero, then consensus 
is determined by the column means of the matrix of \emph{maximum\/} in-forests.

Comparing items\:2 of Theorem~\ref{t_depsop} and Corollary~\ref{c_limInf} we come to the following conclusion: the latent consensus of the protocol with uniform background links coincides with that of the protocol with hub uniformly subordinate to the agents.

\section{Protocols of eigenprojection and orthogonal projection}
\label{s_ortho}

The same purpose, i.e., to develop a concept of latent consensus for the case where classical protocols give no consensus, can be achieved in another way. This can be done by means of the \emph{method of orthogonal projection}~\cite{AgaChe11ARC1,AgaChe15ARC}.
As distinct from the protocols of a dummy hub (Section~\ref{hub_bezv}) and background links (Section~\ref{s_backgr}), where the dependency digraph $\Gamma $ is supplemented by additional links, in the method of orthogonal projection, the link structure is preserved, while the initial state is subject to a minimum necessary correction.
The resulting consensus has the form $\bar{J} S \boldsymbol{x}(0),$ where $S$ is the projection onto the \emph{consensus subspace\/} of the system determined by~$L.$
This consensus generally differs from those of protocol~\eqref{e_modeLK011115c} and the other protocols~\eqref{e_modeLK040715c} studied in this work.

It should be noted that if the dependency digraph contains an in-tree, then all versions of latent consensus coincide with the ordinary
consensus~$\bar{J} \boldsymbol{x}(0).$

The problem of revealing the specific features of two different concepts of latent consensus requires special consideration. In this section, we just provide an example that demonstrates a difference between them.

\begin{example}
For the multi-agent system with Laplacian matrix
\[
L=\left(
  \begin{array}{>{\!\!}r>{\!\!}r>{\!\!}r>{\!\!}r}
    1 &-1 & 0 & 0 \\ [.3em]
   -3 & 3 & 0 & 0 \\ [.3em]
    0 & 0 & 1 &-1 \\[.3em]
    0 & 0 &-4 & 4
  \end{array}
\right),\] we find: 
\begin{gather*}
\bar{J}=\left(
  \begin{array}{>{\!}r>{\!}r>{\!}r>{\!}r}
   0{.}75& 0{.}25&     0&     0\\[.3em]
   0{.}75& 0{.}25&     0&     0\\ [.3em]
        0&      0& 0{.}8& 0{.}2\\[.3em]
        0&      0& 0{.}8& 0{.}2
  \end{array}
\right);\quad
S\approx\left(
  \begin{array}{>{\!\!}r>{\!\!}r>{\!\!}r>{\!\!}r}
    0{.}5690& -0{.}1437&  0{.}4598&  0{.}1149\\ [.3em]
   -0{.}1437&  0{.}9521&  0{.}1533&  0{.}0383\\ [.3em]
    0{.}4598&  0{.}1533&  0{.}5096& -0{.}1226\\ [.3em]
    0{.}1149&  0{.}0383& -0{.}1226&  0{.}9693
  \end{array}
\right);
\\[.1em]
\begin{split}
P_{\mbox{\scriptsize orthog.\,projection}}
=&\,\bar{J} S\approx\bm1\!\cdot\!\big[0{.}3908\;\: 0{.}1303\;\: 0{.}3831\;\: 0{.}0958\big];\\
P_{\mbox{\scriptsize eigenprojection}} =&\,E \bar{J}
=\bm1\!\cdot\!\big[0{.}375\;\: 0{.}125\;\: 0{.}4\;\: 0{.}1\big].
\end{split}\end{gather*}

The two systems of weights obtained for this example do not differ considerably, however, they differently order the weights of the first and the third agents. 
The second order may look more rational, since agent $3$ influences agent $4$ stronger than agent $1$ influences~$2.$ 
On the other hand, a number of examples demonstrate that the method of orthogonal projection gives more ``rewards'' to more homogeneous
components.
\end{example}

\section{Regularization for discrete consensus protocols}
\label{s_discr}

Discrete protocols of consensus seeking can be obtained by replacing differentials with finite differences in the continuous protocols.

\subsection{Protocols with the addition of a hub}
\label{s_discr2}

Consider DeGroot's iterative pooling process~\cite{DeGroot74}:
\begin{gather}\label{e_04072015eq1} 
\boldsymbol{x}^{k+1}=P
\boldsymbol{x}^k,\quad k=0,1,\ldots, 
\end{gather}
where $\boldsymbol{x}^k\in\mathbb{R}^n$ and $P$ is a row-stochastic matrix. 
If $P$ \emph{is not regular\/\footnote{A stochastic matrix is regular iff it is SIA (Stochastic, Indecomposable, Aperiodic).}\/} (i.e., $P$ has at least two eigenvalues of modulus~$1$), then asymptotic consensus is not generally achieved.

As regularized extensions of $P,$ consider the matrices
$$
    Q_{\delta  ,\boldsymbol{v}}=\begin{pmatrix}
        (1-\delta  )P & \delta  \bm1\\[.3em]
         \boldsymbol{v}^{\mathrm T}  &      0
        \end{pmatrix},
$$
where $0<\delta \le 1$ and $\boldsymbol{v}$ is a distribution vector.
$Q_{\delta ,\boldsymbol{v}}$ is regular, since ${v_i>0}$ implies that the $\,i$th column of $Q_{\delta  ,\boldsymbol{v}}^2$
is positive.

Consider the protocols 
\begin{gather}\label{e_deGreG}
\boldsymbol{y}^{k+1}=Q_{\delta
,\boldsymbol{v}}\,\boldsymbol{y}^k,\quad k=0,1,\ldots,
\end{gather} 
where $\boldsymbol{y}^k\in\mathbb{R}^{n+1}.$ These are the counterparts of the continuous protocols \eqref{e_modeLK040715c}. 
Let $\gamma =\delta^{-1}-1.$

\begin{lemma}\label{l_Qsop}
In the above notation\/$,$

$1.$ $(I'-Q_{\delta  ,\boldsymbol{v}})^\vdash
=\bm1'\tfrac1{1+\delta  } \big[\boldsymbol{v}^{\mathrm
T}\big(I+\gamma (I-P)\big)^{-1}\;\:\delta  \big]$$;$

$2.$ $\lim_{\delta  \to+0}(I'-Q_{\delta  ,\boldsymbol{v}})^\vdash
=\bm1' \big[\boldsymbol{v}^{\mathrm T}(I-P)^\vdash\;\:0\big].$
\end{lemma}

Since ${L=I-P}$ is Laplacian, Lemma~\ref{l_Mforest} and the matrix forest theorem are applicable to~$L.$ 
The following Theorem~\ref{t_Qsop} relies on Lemma~\ref{l_Qsop}.

\begin{theorem}\label{t_Qsop}
If $\{\boldsymbol{y}^k\}$ satisfies \eqref{e_deGreG}$,$ then denoting $\,\bar{J}=(I-P)^\vdash$ we have\/$\/:$

$1.$ $\lim\limits_{k\to\infty}\boldsymbol{y}^k =
\bm1'\tfrac1{1+\delta  } \big[\boldsymbol{v}^{\mathrm
T}\big(I+\gamma (I-P)\big)^{-1}\;\;\delta \big]
\boldsymbol{y}^0$\/$;$

$2.$ $\lim\limits_{\delta
\to+0}\lim\limits_{k\to\infty}\boldsymbol{y}^k = \bm1'
\big[\boldsymbol{v}^{\mathrm T}\!\bar{J}\;\;0\big]
\boldsymbol{y}^0$\/$;$

$3.$ If ${\boldsymbol{v}=\tfrac1n\bm1},$ then
${\lim\limits_{\delta
\to+0}\lim\limits_{k\to\infty}\boldsymbol{y}^k
= \bm1'\,\tfrac1n \sum_{j=1}^n\bar{J}_{ \cdot
j}\boldsymbol{y}^0_j},$ where ${\bar{J}_{\cdot j}} 
=\sum_{i=1}^n\bar{J}_{ij}.$
\end{theorem}

In item\:3 of Theorem~\ref{t_Qsop}, we obtain the already familiar to the reader latent consensus: if a hub is subordinate to the agents and infinitesimally influences them, then consensus is determined by averaging the initial agents' states with the weights equal to the column means of~$\bar{J}=(I-P)^\vdash.$

\subsection{Adding background links\/$:$ regularization of PageRank type} \label{s_discr1}

Now for DeGroot's process \eqref{e_04072015eq1}, where $P$ is not regular, consider the protocols 
\begin{gather}\label{e_040715eq3}
\boldsymbol{x}^{k+1}=((1-\delta  )P+\delta   V) \boldsymbol{x}^k,\quad k=0,1,\ldots, 
\end{gather} 
where ${\delta\in\:]0,1]}$, ${V=\bm1\boldsymbol{v}^{\mathrm T}}$, and $\boldsymbol{v}$ is a distribution. 
Such a regularization 
is used\footnote{The iterations of the PageRank method have the form of 
$\big(\bp^{k+1}\big)^{\!{\mathrm T}}=\big(\bp^k\big)^{\!{\mathrm T}}((1-\delta  )P+\delta   V),\; k=0,1,\ldots,$ 
where $\bp^0$ is a distribution. The result is the distribution $\bp^\infty=\lim_{k\to\infty}\bp^k$ (with generally unequal components) not depending on $\bp^0.$ $\big(\bp^\infty\big)^{\!{\mathrm T}}$ coincides with any row of $\lim_{k\to\infty}((1-\delta  )P+\delta V)^k$.} 
in the PageRank method~\cite{BrinPage98,LangvilleMeyer06} (see also~\cite{PolyakTremba12,IshiiTempo14}). The asymptotic behavior of the protocols \eqref{e_040715eq3} is determined by Theorem~\ref{04072015pr1}.

\begin{theorem}\label{04072015pr1}
If a sequence $\{\boldsymbol{x}^k\}$ satisfies~\eqref{e_040715eq3}$,$ then with $\,\bar{J}=(I-P)^\vdash$ it holds that\/$:$

$1.$ $\lim\limits_{k\to\infty}\boldsymbol{x}^k = \bm1
\boldsymbol{v}^{\mathrm T} \big(I+\gamma (I-P)\big)^{-1}
\boldsymbol{x}^0,\,$ where $\:\gamma =\delta  ^{-1}-1$\/$;$

$2.$ $\lim\limits_{\delta
\to+0}\lim\limits_{k\to\infty}\boldsymbol{x}^k = \bm1
\boldsymbol{v}^{\mathrm T}\!\bar{J} \boldsymbol{x}^0$$ ;$

$3.$ If ${\boldsymbol{v}=\tfrac1n\bm1},$ then
${\lim\limits_{\delta
\to+0}\lim\limits_{k\to\infty}\boldsymbol{x}^k =E\bar{J}
\boldsymbol{x}^0 =\bm1\,\tfrac1n \sum_{j=1}^n\bar{J}_{ \cdot
j}\boldsymbol{x}^0_j},$ where ${\bar{J}_{\cdot j}=\sum_{i=1}^n\bar{J}_{ij}.}$
\end{theorem}

Item\:1 of Theorem~\ref{04072015pr1} is equivalent to Eq.~(7.0.1) in~\cite{LangvilleMeyer06}. 

The main conclusion is as follows: item\:3 provides the same latent consensus as the models considered above do.

Finally, one can notice that item\:3 along with Theorem~$2'$ in \cite{AgaChe00'} (interpreting the matrix $\bar{J}=(I-P)^\vdash$) characterizes the limit of iterative regularization for an arbitrary stochastic matrix in~\cite{PolyakTremba12}.

\section{Conclusion}

In this paper, we studied two types of regularization models for the problem of coordination in multi-agent systems. In the models of the first type, the system is supplemented by a dummy agent, ``hub'' (a personalized common ``background''), which uniformly, but very weakly influences the agents and, in turn, depends on them. In the models of the second type, we assume the presence of very weak background links between the initial agents. The paper presents the asymptotic properties of these models. It is demonstrated that in their framework, the concept of latent consensus makes sense. It is the limit consensus as the intensity of the supplemented weak links goes to zero. The latent consensus for the models of both types (wherein both discrete and continuous) turns out to be the same. This sonsensus can be found by averaging the initial states of the agents with the weights equal to the column means of the eigenprojection of the Laplacian matrix $L$ of the dependency digraph (for continuous models) or matrix $I-P,$ where $P$ is the stochastic dependency matrix (for discrete models). An interpretation of these weights is given by Lemma~\ref{l_Mforest} stating the coincidence of the eigenprojection of the Laplacian matrix of a dependency digraph and the normalized matrix of its maximum in-forests.


\appendix{}

\setcounter{lemma}{0}
\def\thelemma{A.\arabic{lemma}}

\PLE{\ref{31102015pr1}}
We need one more lemma.

\begin{lemma}\label{31102015a}
Suppose that $A$ and $C$ satisfy the conditions of Lemma~$\ref{31102015pr1}$. Then for any ${\delta\in\mathbb{R}}$
it holds that ${(A+\delta C)^m =(A+\delta C)(A+\delta I)^{m-1}},$ $\,m=1,2,\ldots.$
\end{lemma}

\PLE{\ref{31102015a}}
For $m=1$ we obtain a trivial identity.
For $m=2$ we have
\begin{gather}\label{e_k2}
(A+\delta C)^2
=A^2+\delta A +\delta C A+\delta^2C 
=A(A+\delta I)+\delta C(A+\delta I)
=(A+\delta C) (A+\delta I),
\end{gather} 
i.e., the desired statement holds true.
Let it be true for ${m=k>1}$. Now we prove it for 
${m=k+1}$. Due to the induction hypothesis and \eqref{e_k2}, we obtain
\[
 (A+\delta C)^{k+1}
=(A+\delta C)(A+\delta C)(A+\delta I)^{k-1}
=(A+\delta C)(A+\delta I)^k.
\] Lemma~\ref{31102015a} is proved.

Now we use the series representation of the matrix exponent:
\begin{gather*}
 \begin{split}                   e^{-(A+\delta C)t}
                    &=\sum_{k=0}^\infty\frac{(-t)^k(A+\delta C)^k}{k!}\\
                    &= I-\frac{t(A+\delta C)}{1!}+\frac{t^2(A+\delta C)^2}{2!}
                     - \frac{t^3(A+\delta C)^3}{3!}+\ldots
                     \end{split}
\end{gather*}

By Lemma~\ref{31102015a} we have
\begin{eqnarray}\label{e_31102015}
\quad
e^{-(A+\delta C)t}
=I-(A+\delta C)\left(\frac{tI}{1!}-\frac{t^2(A+\delta I)}{2!}+\frac{t^3(A+\delta I)^2}{3!}+\ldots\right)\!.
\end{eqnarray}

Since by the conditions of Lemma~\ref{31102015pr1} $A+\delta I$ is invertible, \eqref{e_31102015} can be rewritten in the form {\arraycolsep=.5pt
\begin{gather}\label{e_31102015-3}
\begin{gathered}
\begin{split}
e^{-(A+\delta C)t}
&=I+(A+\delta C)(A+\delta I)^{-1}
\left(-I+I-\frac{t(A+\delta I)}{1!}+\frac{t^2(A+\delta I)^2}{2!}-\ldots\!\right)\\
&= I+(A+\delta C)(A+\delta I)^{-1}\left(e^{-(A+\delta I)t}-I\right).
\end{split}
\end{gathered}
\end{gather}

Since the real parts of the eigenvalues of $A+\delta I$ are positive, $\lim_{t\to\infty}e^{-(A+\delta I)t}=0$ holds. From
\eqref{e_31102015-3} we obtain
\begin{gather*}\begin{split}
\lim_{t\to\infty}
e^{{\hspace*{1pt}-\hspace*{1pt}}(A{\hspace*{1pt}+\hspace*{1pt}}\delta
C) t}
 &=                I{\hspace*{1pt}-\hspace*{1pt}}(A{\hspace*{1pt}+\hspace*{1pt}}\delta C)(A{\hspace*{1pt}+\hspace*{1pt}}\delta I)^{{\hspace*{1pt}-\hspace*{1pt}}1}
  =       I{\hspace*{1pt}-\hspace*{1pt}}(A{\hspace*{1pt}+\hspace*{1pt}}\delta I{\hspace*{1pt}+\hspace*{1pt}}\delta C{\hspace*{1pt}-\hspace*{1pt}}\delta I)(A{\hspace*{1pt}+\hspace*{1pt}}\delta I)^{{\hspace*{1pt}-\hspace*{1pt}}1}\\
 &=              I{\hspace*{1pt}-\hspace*{1pt}}I{\hspace*{1pt}+\hspace*{1pt}}\delta (I{\hspace*{1pt}-\hspace*{1pt}}C)(A{\hspace*{1pt}+\hspace*{1pt}}\delta I)^{{\hspace*{1pt}-\hspace*{1pt}}1}
  =      (I{\hspace*{1pt}-\hspace*{1pt}}C)\big(I{\hspace*{1pt}+\hspace*{1pt}}\delta^{{\hspace*{1pt}-\hspace*{1pt}}1}  A\big)^{{\hspace*{1pt}-\hspace*{1pt}}1}.\end{split}
\end{gather*}
Lemma~\ref{31102015pr1} is proved.

\PLE{\ref{l_sopL}}
Since $A$ is Laplacian, its spectrum  lies in the right half-plane bounded by the vertical axis \cite[proposition\:9]{AgaChe01AiT}, therefore, the real parts of the eigenvalues of $\delta I+A$ are positive. Thereby the condition of Lemma~\ref{31102015pr1} are fulfilled, thus, so does its conclusion. Comparing it with that of Lemma~\ref{l_sopLe} (whose proof is given below) and having in mind that
$A+\delta C$ is a Laplacian matrix, we obtain the desired statement.
Lemma~\ref{l_sopL} is proved.

\PLE{\ref{l_sopLe}}
All the solutions of the equation 
${\dot {\boldsymbol{x}}(t)=-L\,\boldsymbol{x}(t)}$ 
satisfy the identity \cite[Eq.\,(43) of Chapter~5]{Gantmacher59} 
${\boldsymbol{x}(t)=e^{-L t} \boldsymbol{x}(0)} $, 
from which in the case of existence of the limit at ${t\to\infty},$ it follows that
${\lim_{t\to\infty}\boldsymbol{x}(t)=\lim_{t\to\infty}e^{-L t} \boldsymbol{x}(0)}$. 
Comparing the last expression with Proposition~\ref{p_limInf} and taking into account the arbitrariness of
$\boldsymbol{x}(0)$, one has $L^\vdash = \lim_{t\to\infty}e^{-L t}.$
Lemma~\ref{l_sopLe} is proved.

\PLE{\ref{02072015b}} 
Set ${A\!:=L_0+H_{\boldsymbol{v}}}$, ${C\!:=H_I}$. It can be easily checked that ${AC=A}$ and ${C^2=C}$. By Lemma~\ref{l_sopL},
${(L_0+H_{\delta ,\boldsymbol{v}})^\vdash 
= (I'-H_I)(I'+\tfrac1\delta L_0+\tfrac1\delta H_{\boldsymbol{v}})^{-1}}$. 
Since ${I'+\tfrac1\delta L_0+\tfrac1\delta H_{\boldsymbol{v}} 
= \big(I'+\tfrac1\delta L_0\big)\big(I'+\tfrac1\delta H_{\boldsymbol{v}}\big)}$, we have
\begin{gather*}
\begin{split}
(L_0+H_{\delta  ,\boldsymbol{v}})^\vdash
&= (I'-H_I)\left(I'+\dfrac1\delta   H_{\boldsymbol{v}}\right)^{-1}\left(I'+\dfrac1\delta   L_0\right)^{-1}\\
&= \big[0_{(n+1)\times n}\;\,\bm1'\big]\begin{pmatrix}I&\bm0\\[.6em]
                                                \dfrac1{s+\delta  }\boldsymbol{v}^{\mathrm T}&\;\dfrac\delta  {s+\delta  }\end{pmatrix}\!
                                        \begin{pmatrix}
                                        \left(I+\dfrac1\delta   L\right)^{-1}&\bm0\\[1em]
                                        \bm0^{\mathrm T}         &   1
                                        \end{pmatrix}\\
&= \bm1'\dfrac1{s+\delta  }\left[\boldsymbol{v}^{\mathrm
T}\left(I+\dfrac1\delta   L\right)^{-1}\;\;\delta  \right].
\end{split}
\end{gather*}
Lemma~\ref{02072015b} is proved.

\PLE{\ref{l_SopNs}}
From Lemmas~\ref{02072015b} and~\ref{02072015th1} one derives:
\begin{gather*}
\begin{split}
\lim_{\delta  \to+0}(L_0+H_{\delta  ,\,\boldsymbol{v}})^\vdash 
&=\lim_{\delta  \to+0}\bm1'\frac1{s+\delta  }
\left[\boldsymbol{v}^{\mathrm T}\left(I+\frac1\delta
L\right)^{-1}\;\;\, \delta  \right]\\
&= \lim_{\delta
\to+0}\bm1'\frac1{1+\delta  /s}
\left[\frac1s\boldsymbol{v}^{\mathrm T}\bar{J}\;\;\,\frac\delta
s\right]
= \bm1'\Big[\tilde{\boldsymbol{v}}^{\mathrm T}\!\bar{J}\;\;\, 0\Big].
\end{split}
\end{gather*}
Lemma~\ref{l_SopNs} is proved.

\PLE{\ref{l_PRsop}}
1. For ${A=L}$ and ${C=D},$ the conditions of Lemma~\ref{l_sopL} are fulfilled, therefore, 
${(L+\delta   D)^\vdash 
=(I-D)\big(I+\tfrac1\delta L\big)^{-1} 
= \bm1 \boldsymbol{v}^{\mathrm T} \big(I+\tfrac1\delta L\big)^{-1}}$.

2. Using Lemma~\ref{02072015th1} we obtain 
${\lim\limits_{\delta \to0}(L+\delta   D)^\vdash 
=\lim\limits_{\delta \to0}\bm1 \boldsymbol{v}^{\mathrm T} \big(I+\tfrac1\delta L\big)^{-1} 
=\bm1 \boldsymbol{v}^{\mathrm T}\!\bar{J}}
$. 
Lemma~\ref{l_PRsop} is proved.

\PLE{\ref{l_Qsop}}
Set
$$
    P_0=\begin{pmatrix}
          P     & \bm0\\[.3em]
         \bm0^{\mathrm T}& 0
        \end{pmatrix},\;\;
    P_1=\begin{pmatrix}
        0_{n\times n} & \delta  \bm1\\[.3em]
        \bm0^{\mathrm T}       &      0
        \end{pmatrix},\;\;
P_2=  \begin{pmatrix}
       0_{n\times n} & \bm0\\[.3em]
               \boldsymbol{v}^{\mathrm T} &    0
        \end{pmatrix},\;\mbox{ and }\;
\hat I=I'-I_0.
$$

Then ${Q_{\delta  ,\boldsymbol{v}}=(1-\delta  )P_0+\delta P_1+P_2}$. Using \eqref{e_W1}, \eqref{e_modeLK040715c}, and the 
notation ${L_0\!:=(1-\delta )(I_0-P_0)}$ ($L_0$ is a Laplacian matrix!), we obtain
\begin{gather*}
\begin{split}I'- Q_{\delta  ,\boldsymbol{v}}
&= (1-\delta  )(I_0-P_0)+\delta  (I_0-P_1)+(\hat I-P_2)\\
&= (1-\delta  )(I_0-P_0)+\delta   H_I     +H_{\boldsymbol{v}} = L_0
+\delta   H_I     +H_{\boldsymbol{v}} =              L_0 +H_{\delta
,\boldsymbol{v}}.\end{split}
\end{gather*}
Now the desired statement 
is derived from Lemma~\ref{02072015b} (with substitution ${L=(1-\delta  )(I-P)}$) and
Lemma~\ref{02072015th1}.
Lemma~\ref{l_Qsop} is proved.

\PTH{\ref{t_Qsop}}
Since $Q_{\delta  ,\boldsymbol{v}}$ is a regular stochastic matrix, the limit of its powers coincides \cite{Meyer75,Rothblum76ai} with the eigenprojection of ${I'-Q_{\delta,\boldsymbol{v}}}$. Thus,
${\lim_{k\to\infty}\boldsymbol{y}^k 
= \lim_{k\to\infty}Q_{\delta ,\boldsymbol{v}}^k \boldsymbol{y}^0
= (I'-Q_{\delta ,\boldsymbol{v}})^\vdash\boldsymbol{y}^0}$, and the theorem follows from Lemma~\ref{l_Qsop}.
Theorem~\ref{t_Qsop} is proved.

\PTH{\ref{04072015pr1}}
1. Since the stochastic matrix ${(1-\delta )P+\delta V}$ is regular, the limit of its powers coincides \cite{Meyer75,Rothblum76ai} with
the eigenprojection of 
${ I-(1-\delta  )P-\delta V =(1-\delta  )(I-P)+\delta  (I-V)
=L+\delta  D} $, where
${L\!:=(1-\delta )(I-P)}$ is a Laplacian matrix and ${D=I-V}$ (see~\eqref{e_modeL-PR}).
Therefore, by Lemma~\ref{l_PRsop} the required limit of the powers equals ${E\big(I+\tfrac1\delta   L\big)^{-1}}$, which implies item\:1 of Theorem~\ref{04072015pr1}.

2. The desired statement (along with item\:3) also follows from Lemma~\ref{l_PRsop} in view of the fact that the eigenprojection is preserved under the multiplication of a matrix by a nonzero scalar.
Theorem~\ref{04072015pr1} is proved.



\begin{thebibliography}{99}
 \itemsep0mm
\parsep0mm

\bibitem{Olfati-SaberMurray04} 
  {Olfati-Saber,~R. and Murray,~R.M.},
  {Consensus Problems in Networks of Agents with Switching Topology and Time-delays},
  \emph{IEEE Trans. Automat. Control},
  2004,
  vol.~49, no.~9,
  pp.~1520--1533.
  
\bibitem{MesbahiEgerstedt10book} 
  {Mesbahi,~M. and Egerstedt,~M.},
  \emph{Graph Theoretic Methods in Multiagent Networks}.
  Princeton: Princeton University Press, 2010.

\bibitem{AgaChe11ARC1}
  {Agaev,~R.P. and Chebotarev,~P.Yu.},
  {The Projection Method for Reaching Consensus and the Regularized Power Limit of a Stochastic Matrix},
  \emph{Automat. Remote Control},
  2011,
  vol.~72, no.~12,
  pp.~2458--2476.

\bibitem{AgaChe15ARC}
  {Agaev,~R.P. and Chebotarev,~P.Yu.},
  {The Projection Method for Continuous-Time Consensus Seeking}, 
  \emph{Automat. Remote Control}, 
  2015,
  vol.~76, no.~8, 
  pp.~1436--1445.

\bibitem{Lin99} 
    {Lin, Z.}, 
    \emph{Low Gain Feedback}. 
    London: Springer, 1999.

\bibitem{SeoShimBack09}
{Seo,~J.H., Shim,~H. and Back,~J.},
{Consensus of High-Order Linear Systems Using Dynamic Output Feedback Compensator: Low Gain Approach}, 
\emph{Automatica}, 
2009. 
vol.~45, no.~11, 
pp.~2659--2664.

\bibitem{SuChenLamLin13} 
{Su,~H., Chen.~M.Z., Lam,~J. and Lin,~Z.}, 
{Semi-Global Leader-Following Consensus of Linear Multi-Agent Systems with Input Saturation via Low Gain Feedback}, 
\emph{IEEE Trans. Circuits and Systems I: Regular Papers}, 
2013, 
vol.~60, no.~7, 
pp.~1881--1889.

\bibitem{RenCao11distributed} 
{Ren,~W. and Cao~Y.} 
\emph{Distributed Coordination of Multi-Agent Networks: Emergent Problems, Models, and Issues}. 
London: Springer, 2011.

\bibitem{CheAga02ap} 
  {Chebotarev,~P. and Agaev,~R.},
  {Forest Matrices around the {Laplacian} Matrix},
  \emph{Linear Algebra and its Applications},
  2002,
  vol.~356,
  pp.~253--274.
  
\bibitem{CheAga14IEEETAC} 
  {Chebotarev,~P. and Agaev,~R.},
  {The Forest Consensus Theorem},
  \emph{IEEE Trans. Automat. Control},
  2014,
  vol.~59, no.~9,
  pp.~2475--2479.

\bibitem{AgaChe00'}
  {Agaev,~R.P. and Chebotarev,~P.Yu.},
  {The Matrix of Maximum {Out} Forests of a Digraph and its Applications},
  \emph{Automat. Remote Control},
  2000,
  vol.~61, no.~9,
  pp.~1424--1450.

\bibitem{CheSha97'}
  {Chebotarev,~P.Yu. and Shamis,~E.V.}, 
  {The Matrix-Forest Theorem and Measuring Relations in Small Social Groups}, 
  \emph{Automat. Remote Control},
  1997, 
  vol.~58, no.~9, 
  pp.~1505--1514.

%
\bibitem{DeGroot74} 
  {DeGroot,~M.H.},
  {Reaching a Consensus},
  \emph{J. Amer. Statistical Association},
  1974,
  vol.~69, no.~345,
  pp.~118--121.

\bibitem{BrinPage98} 
    {Brin,~S. and Page,~L.}, 
    {The Anatomy of a Large-Scale Hypertextual Web Search Engine},
    \emph{Computer Networks and ISDN Systems}, 
    1998, 
    vol.~30, 
    pp.~107--117.

\bibitem{LangvilleMeyer06} 
    {Langville,~A.N. and Meyer,~C.D.},
    \emph{Google's PageRank and Beyond: The Science of Search Engine Rankings}.
    Princeton: Princeton Univ. Press, 2006.

\bibitem{PolyakTremba12}
    {Polyak,~B.T. and Tremba,~A.A.},
    {Regularization-Based Solution of the PageRank Problem for Large Matrices}, 
    \emph{Automat. Remote Control},
    2012, 
    vol.~73, no.~11.,
    pp.~1877--1894.

\bibitem{IshiiTempo14} 
    {Ishii,~H. and Tempo,~R.}, 
    {The PageRank Problem, Multiagent Consensus, and Web Aggregation: A Systems and Control Viewpoint}, 
    \emph{IEEE Control Syst. Magazine}, 
    2014,
    vol.~34, no.~3, 
    pp.~34--53.

\bibitem{AgaChe01AiT}
  {Agaev,~R.P. and Chebotarev,~P.Yu.},
   {Spanning Forests of a Digraph and Their Applications},
  \emph{Automat. Remote Control},
    2001,
    vol.~62, no.~3, 
    pp.~443--466.

\bibitem{Gantmacher59}
  {Gantmacher,~F.R.},
  \emph{The Theory of Matrices}.\hskip 1em plus 0.5em minus
  0.4em\relax New York: Chelsea, 1959.

%

\bibitem{Meyer75}
{Meyer,~C.D.,~Jr.}, 
{The Role of the Group Generalized Inverse in the Theory of Finite {Markov} Chains},
\emph{SIAM Review}, 
1975,
vol.~17, no.~3, 
pp.~443--464.

\bibitem{Rothblum76ai} 
{Rothblum,~U.G.}, 
{Computation of the Eigenprojection of a Nonnegative Matrix at Its Spectral Radius},     
\emph{Stochastic Systems: Modeling, Identification and Optimization~II, ser. Mathematical Programming Study}, R.~J.-B.~Wets, Ed. 
Amsterdam: North-Holland, 1976, vol.~6, pp.~188--201.
\end{thebibliography}
\end{document}